\documentclass[a4paper,11pt]{article}
\usepackage[top=2.5cm,bottom=2.5cm,left=2.2cm,right=2.2cm]{geometry}
\usepackage{amsfonts}
\usepackage{mathrsfs,amscd,amssymb,amsthm,amsmath,bm,graphicx,psfrag,subfigure,url,mathtools}
\usepackage{pict2e}
\usepackage{stfloats}
\usepackage{psfrag,amsmath}
\usepackage{tikz}
\usepackage{indentfirst}
\usepackage{hyperref}
\usepackage{bookmark}
\usepackage{enumerate}
\usepackage{latexsym,euscript,epic,eepic,color}
\usepackage{multirow}
\usepackage{multicol}
\usepackage{longtable}
\usepackage{adjustbox}
\usepackage[all]{xy}
\usepackage{setspace}
\usepackage{epstopdf}
\allowdisplaybreaks
\usepackage{authblk}
\usepackage{pifont}

\makeatletter

\renewcommand{\@seccntformat}[1]{{\csname the#1\endcsname}{\normalsize .}\hspace{.5em}}
\makeatother

\usepackage{ifpdf}
\usepackage{indentfirst}

\def \[{\begin{equation}}
\def \]{\end{equation}}

\newcommand{\ex}{{\rm ex}}
\newcommand{\Ex}{{\rm Ex}}
\newtheorem{thm}{Theorem}[section]

\newtheorem{claim}{Claim}

\newtheorem{lem}[thm]{Lemma}

\newtheorem{pb}{Problem}

\newenvironment{wst}
{\setlength{\leftmargini}{1.5\parindent}
 \begin{itemize}
 \setlength{\itemsep}{-1.1mm}}
{\end{itemize}}

\begin{document}
\baselineskip=0.23in

\title{\bf The exact Tur\'an number of generalized book graph $B_{r,k}$ in non-$r$-partite graphs%Refinement on Tur\'an number of the generalized book graph
\thanks{S.L. financially supported by the National Natural Science Foundation of China (Grant Nos. 12171190, 11671164), the Special Fund for Basic Scientific Research of Central Colleges (Grant Nos. CCNU25JC006, CCNU25JCPT031) and the Open Research Fund of Key Laboratory of Nonlinear Analysis \& Applications (CCNU), Ministry of Education of China (Grant No. NAA2025ORG010)\\[3pt] \hspace*{5mm}{\it Email addresses}:
ytyumath@sina.com (Y. Yu),\ li@ccnu.edu.cn (S. Li)}}

\author[1]{Yuantian Yu}% Author
\author[,2,3]{Shuchao Li\thanks{Corresponding author}}
%Communication author
\affil[1]{School of Science, East China University of Technology, Nanchang 330013, China}
\affil[2]{School of Mathematics and Statistics, and Hubei Key Lab--Math. Sci.,\linebreak Central China Normal University, Wuhan 430079, China}% bit
\affil[3]{Key Laboratory of Nonlinear Analysis \& Applications (Ministry of Education),\linebreak Central China Normal University, Wuhan 430079, China}
\date{\today}
\maketitle

\begin{abstract}
Given a graph $H,$ we say that a graph is \textit{$H$-free} if it does not contain $H$ as a subgraph. The Tur\'an number $\ex(n,H)$ of $H$ is the maximum number of edges in an $n$-vertex $H$-free graph, the set of all the corresponding extremal graphs is denoted by $\Ex(n, H)$. The study of Tur\'an number of graphs is a central topic in extremal graph theory. A graph is \textit{color-critical} if it contains an edge whose deletion reduces its chromatic number. Simonovits showed that if $H$ is a color-critical graph of chromatic number $r+1,$ then for sufficiently large $n,$ $\Ex(n, H)=\{T_r(n)\},$ the $r$-partite Tur\'an graph of order $n.$ Given a color-critical graph $H$ with chromatic number $r+1,$ it is interesting to determine $H$-free non-$r$-partite graphs with maximum number of edges. For a graph $H$ with chromatic number $r+1,$ denote $\ex_{r+1}(n,H)$ the maximum number of edges in non-$r$-partite $H$-free graphs of order $n,$ the set of all non-$r$-partite $H$-free graphs of order $n$ and size $\ex_{r+1}(n,H)$ is denoted by $\Ex_{r+1}(n, H)$. For $r\geq 3,\,k\geq1,$ the generalized book graph \({B}_{r,k}\) is a graph obtained by joining every vertex of $K_r$ to every vertex of an independent set of size \(k\). Note that \({B}_{r,k}\) is a color-critical graph of chromatic number $r+1.$ In this paper, based on the stability theory and local structure characterization, the exact value of $\ex_{r+1}(n,B_{r,k})$ is determined and all the corresponding extremal graphs are identified, where $r\geq 3,\,k\geq1$ and $n$ is sufficiently large. 
\vspace{2mm}

\end{abstract}

%%%%%%%%%%%%%

\section{\normalsize Introduction}\label{s1}
In this paper, we consider only simple, undirected, and finite graphs. Unless otherwise stated, we follow the traditional notation and terminology (see, for instance, Bollob\'as \cite{BB1998}, West \cite{DW1996}).

For a graph $G=(V(G),E(G)),$ we use $|V(G)|$ and $e(G):=|E(G)|$ to denote the \textit{order} and the \textit{size} of $G,$ respectively. With no confusion, we also use the \textit{size} to denote the cardinality of a set. As usual, let $K_n$ and $C_n$ be the complete graph and cycle of order $n$, respectively. 
A simple complete $r$-partite graph on $n$ vertices whose parts are of sizes $t_1,t_2,\ldots, t_r$ is denoted by $K_{t_1,t_2,\ldots,t_r},$ in which $t_1+\cdots+t_r=n$.

Given a graph $H,$ we say that a graph is \textit{$H$-free} if it does not contain $H$ as a subgraph. The Tur\'{a}n type problem asks: determine the maximum number of edges, $\ex(n,H)$, of an $n$-vertex $H$-free graph; characterize the set, $\Ex(n,H)$, of all the $n$-vertex $H$-free graph of size $\ex(n,H)$. The value of $\ex(n,H)$ is called the \textit{Tur\'an number} of $H$, the entry in $\Ex(n,H)$ is called the \textit{extremal graph} for $H$. The research on the Tur\'an-type problem attracts much attention, and it has become to be one of the most attractive fundamental problems in extremal graph theory (see \cite{FS2013,V2011} for surveys). 

Let $T_r(n)$ denote the complete $r$-partite graph of order $n$ whose parts are of equal or almost
equal. That is, $T_r(n)=K_{t_1,t_2,\ldots,t_r}$ with $\sum_{i=1}^r t_i=n$ and $|t_i-t_j|\leq 1$ for $i\neq j.$ In 1941, Tur\'an \cite{T1941} determined the Tur\'an number of $K_{r+1},$ which extended the result of Mantel \cite{M1907}. 
\begin{thm}[\cite{T1941}]\label{thm1.1}
For positive integers $r,n,$ we have $\ex(n,K_{r+1})=e(T_r(n))$ and $\Ex(n,K_{r+1})=\{T_r(n)\}.$ 
\end{thm}

A graph is said to be \textit{properly coloured} if each vertex is coloured so that the end vertices of each edge have different colours. The chromatic number $\chi(G)$ is the minimum $k$ such that $G$ can be properly coloured by $k$ colours. A graph $G$ is color-critical if $G$ contains an edge $e$ such that $\chi(G - e) < \chi(G).$ Note that $K_{r+1}$ is color-critical with chromatic number $r+1,$ Simonovits extended Theorem \ref{thm1.1} to color-critical graphs. 
\begin{thm}[\cite{Simo1968,Simo1974}]\label{thm1.2}
For a positive integer $r$ and a graph $H,$ if $H$ is color-critical with $\chi(H)=r+1$, then for sufficiently large $n,$ we have $\ex(n,H)=e(T_r(n))$ and $\Ex(n,H)=\{T_r(n)\}.$ 
\end{thm}

Let $H$ be a graph with $\chi(H)=r+1,$ denote $\ex_{r+1}(n,H)$ the maximum number of edges in non-$r$-partite $H$-free graphs of order $n,$ and let $\Ex_{r+1}(n, H)$ be the set of all non-$r$-partite $H$-free graphs of order $n$ and size $\ex_{r+1}(n,H)$. 

By Theorem \ref{thm1.2}, given a color-critical graph $H$ with $\chi(H)=r+1$, for sufficiently large $n,$ the unique extremal graph for $H$ is $T_r(n),$ which is $r$-partite. It is interesting to consider the following problem.
\begin{pb}\label{pb1}
Given a color-critical graph $H$ with $\chi(H)=r+1$. Determine the value of $\ex_{r+1}(n,H)$ and characterize all the extremal graphs in $\Ex_{r+1}(n, H).$
\end{pb}

Amin et al. \cite{AFGS}, Caccetta and Jia \cite{CJ2002},  Erd\H{o}s \cite{E1955}, Kang and Pikhurko \cite{KP2005}, independently, solved Problem \ref{pb1} for $C_3.$ For $k\geq 2,$ Ren et al. \cite{RWWY} solved Problem \ref{pb1} for $C_{2k+1}.$   

Let \({v}_{1}{v}_{2}{v}_{3}{v}_{4}{v}_{5}{v}_{1}\) be a 5-cycle, we denote by \({C}_{5}\left\lbrack  {{n}_{1},{n}_{2},{n}_{3},{n}_{4},{n}_{5}}\right\rbrack\) a blow-up of \({C}_{5}\), which is obtained from \({v}_{1}{v}_{2}{v}_{3}{v}_{4}{v}_{5}{v}_{1}\) by the following way: for each \(1\leq i\leq 5\), we replace \({v}_{i}\) with an independent set \({I}_{i}\) of size \({n}_{i}\); for each \(1\leq i \leq 4\), we add all possible edges between \({I}_{i}\) and \({I}_{i + 1}\), and all possible edges between \({I}_{1}\) and \({I}_{5}\). For convenience, define
$$
\mathcal{C}_5^\ast=\left\{{C}_{5}\left\lbrack  {{n}_{1},{n}_{2},{n}_{3},{n}_{4},{n}_{5}}\right\rbrack: {n}_{1}+{n}_{2}+{n}_{3}+{n}_{4}+{n}_{5}\geq 5\right\}.
$$
For even $n\geq 6,$ define
\begin{align*}
\mathcal{C}_5^1[n]&=\left\{{C}_{5}\left\lbrack  {\frac{n}{2}-2,t,1,1,\frac{n}{2}-t}\right\rbrack: 1\leq t\leq \frac{n}{2}-1\right\}, \\
\mathcal{C}_5^2[n]&=\left\{{C}_{5}\left\lbrack  {\frac{n}{2}-1,t,1,1,\frac{n}{2}-t-1}\right\rbrack: 1\leq t\leq \frac{n}{2}-2\right\}.
\end{align*}
For odd $n\geq 5,$ define
$$
\mathcal{C}_5^3[n]=\left\{{C}_{5}\left\lbrack  {\frac{n-1}{2}-1,t,1,1,\frac{n-1}{2}-t}\right\rbrack: 1\leq t\leq \frac{n-1}{2}-1\right\}.
$$

%\begin{thm}[\cite{AFGS,KP2005}]\label{thm1.3}
%Let $n,q,p$ be integers satisfy $n=2q+p\geq5,\,p=0,1.$ Then 
%\begin{wst}
%\item[{\rm (a)}]$\ex_{3}(n,K_3)=\frac{n^2}{4}-\frac{n}{2}+\frac{p}{4}+1;$ 
%\item[{\rm (b)}]$\Ex_{3}(n,K_3)=\left\{
%                                 \begin{array}{ll}
%                                    \mathcal{C}_5^1[n]\cup \mathcal{C}_5^2[n], & \text{if $n$ is even}; \\
%                                                 \mathcal{C}_5^3[n], & \text{if $n$ is even}.
%                                               \end{array}
%                                             \right.
%$
%\end{wst}
%\end{thm}
Let $H_1$ and $H_2$ be two graphs, define $H_1\cup H_2$ to be their disjoint union, and $H_1\vee H_2$ to be their \textit{join}, which is obtained from $H_1\cup H_2$ by connecting each vertex of $H_1$ with each vertex of $H_2$ by an edge. Let $n,r,q,p$ be integers satisfy $r\geq3,\,n=rq+p\geq r+3,\,0\leq p\leq r-1.$ Then define 
\begin{align*}
\mathcal{G}_1[n]&=\left\{F\vee T_{r-2}(q(r-2)+p+1):F\in \mathcal{C}_5^3[2q-1]\right\}; \\
\mathcal{G}_2[n]&=\left\{F\vee T_{r-2}(q(r-2)+p):F\in \mathcal{C}_5^1[2q]\cup\mathcal{C}_5^2[2q]\right\}; \\
%\mathcal{G}_3&=\left\{F\vee T_{r-2}(q(r-2)+p):F\in \mathcal{C}_5^2[2q]\right\}; \\
\mathcal{G}_3[n]&=\left\{F\vee T_{r-2}(q(r-2)+p-1):F\in \mathcal{C}_5^3[2q+1]\right\}.
\end{align*}

For $r\geq3,$ Amin et al. \cite{AFGS} and Kang, Pikhurko \cite{KP2005} solved Problem \ref{pb1} for $K_{r+1},$ independently.
\begin{thm}[\cite{AFGS,KP2005}]\label{thm1.4}
Let $n,r,q,p$ be integers satisfy $r\geq3,\,n=rq+p\geq r+3,\,0\leq p\leq r-1.$ Then 
\begin{wst}
\item[{\rm (a)}]$\ex_{r+1}(n,K_{r+1})=\left( {1 - \frac{1}{r}}\right)  \cdot  \frac{{n}^{2}}{2} - \frac{n}{r} + \frac{p\left( {p + 2}\right) }{2r} - \frac{p}{2} + 1;$ 
\item[{\rm (b)}]$\Ex_{r+1}(n,K_{r+1})=\left\{
                                        \begin{array}{ll}
                            \left\{C_5\vee T_{r-2}(n-5)\right\}, & \text{if $q=1,2;$} \\
                       \mathcal{G}_1[n]\cup\mathcal{G}_2[n], &  \text{if $q\geq 3$ and $p=0;$} \\
\mathcal{G}_1[n]\cup\mathcal{G}_2[n]\cup\mathcal{G}_3[n], &  \text{if $q\geq 3$ and $1\leq p\leq r-3;$} \\
\mathcal{G}_2[n]\cup\mathcal{G}_3[n], &  \text{if $q\geq 3$ and $p=r-2;$} \\
\mathcal{G}_3[n], &  \text{if $q\geq 3$ and $p=r-1.$}
                                        \end{array}
                                      \right.
$
\end{wst}
\end{thm}

For positive integers $k,r$ with $r\geq3,$ define \({B}_{r,k}={K}_{r}\vee kK_1\) to be the generalized book graph. Clearly, \({B}_{r,k}\) is a color-critical graph with $\chi({B}_{r,k})=r+1$. In this paper, we consider Problem \ref{pb1} for generalized book graphs. 
The main result of this paper is presented in the following.
\begin{thm}\label{thm1}
Let \(r,k,n\) be integers, where \(r \geq  3,\,k \geq  1\) and \(n\) is sufficiently large. If \(n = {qr} + p\) for some integers \(p,\,q\) with \(0 \leq  p \leq  r - 1\), %and \(G\) is a non-\(r\)-partite \({B}_{r,k}\) -free graph of order \(n\), 
then
\begin{wst}
\item[{\rm (a)}]$\ex_{r+1}(n,{B}_{r,k})=\left( {1 - \frac{1}{r}}\right)  \cdot  \frac{{n}^{2}}{2} - \frac{n}{r} + \frac{p\left( {p + 2}\right) }{2r} - \frac{p}{2} + 1;$ 
\item[{\rm (b)}]$\Ex_{r+1}(n,{B}_{r,k})=\left\{
                                        \begin{array}{ll}
                       \mathcal{G}_1[n]\cup\mathcal{G}_2[n], &  \text{if $p=0;$} \\
\mathcal{G}_1[n]\cup\mathcal{G}_2[n]\cup\mathcal{G}_3[n], &  \text{if $1\leq p\leq r-3;$} \\
\mathcal{G}_2[n]\cup\mathcal{G}_3[n], &  \text{if $p=r-2;$} \\
\mathcal{G}_3[n], &  \text{if $p=r-1.$}
                                        \end{array}
                                      \right.
$
\end{wst}
\end{thm}
%\begin{thm}
%Let \(r \geq  3,k \geq  1\) and \(n\) be sufficiently large. If \(G\) is a non-1-partite \({B}_{r,k}\) -free graph of order \(n\) , then \(P\left( G\right)  \leq  P\left( {{Y}_{r}\left( {n,1}\right) }\right)\) , with equality if and only if \(G \cong  {Y}_{r}\left( {n,1}\right)\) .
%\end{thm}
\noindent{\bf Outline of our paper.}\ \ In the remainder of this section, we introduce some necessary notations and terminologies. In Section \ref{s2}, we give some necessary preliminaries. In Section \ref{s3}, by applying the method of stability, we progressively refine the structure of our extremal graphs and complete the proof of Theorem~\ref{thm1}. Some concluding remarks are given in the last section.\\

\noindent{\bf Notations and terminologies.}\ \ Let $G$ be a graph, we say that two vertices $u$ and $v$ in $G$ are \textit{adjacent} (or \textit{neighbours}) if they are joined by an edge. If $uv\in E(G),$ then let $G-uv$ (resp. $G-u$) denote the graph obtained from $G$ by deleting edge $uv$ (resp. vertex $u$) (this notation is naturally extended if more than one edge (resp. vertex) is deleted). Similarly, if $uv\notin E(G),$ then let $G+uv$ denote the graph obtained from $G$ by adding an edge joining $u$ and $v.$ For two disjoint vertex subsets $V_1$ and $V_2$ of $V(G),$ denote by $G[V_1]$ a subgraph of $G$ induced on $V_1$ and $G[V_1,V_2]$ a subgraph of $G$ induced on the edges between $V_1$ and $V_2$. Then the number of edges of $G[V_1]$ and $G[V_1,V_2]$ can be abbreviated to $e(V_1)$ and $e(V_1,V_2)$, respectively. The \textit{set of neighbours} of a vertex $v$ (in $G$) is denoted by $N_G(v)$, its size is the \textit{degree} of $v$ (in $G$), denoted by $d_G(v).$ Further, for a vertex $v\in V(G)$ and a subset $W\subseteq V(G),$ denote $N_W(v)=N_G(v)\cap W,\,d_W(v)=|N_W(v)|.$ For an positive integer $t,$ the set $\{1,2,\ldots,t\}$ is abbreviated as $[t].$ 

\section{\normalsize Preliminaries}\label{s2}
By a simple calculation, we have the following fact.
\begin{lem}\label{lem2.1}  
For $n\geq r\geq 2,$ if $G$ is an $r$-partite graph of order $n,$ then $e(G)\leq e\left(T_r(n)\right).$ Further, the size of $T_r(n)$ satisfies:
$$
\left(1-\frac{1}{r}\right)\frac{n^2}{2}-\frac{r}{8}\leq e\left(T_r(n)\right)
\leq \left(1-\frac{1}{r}\right)\frac{n^2}{2}.
$$ 
\end{lem}
\begin{lem}[\cite{CFTZ}]\label{lem2.2}
Let $V_1,V_2,\ldots,V_t$ be $t$ finite sets. Then
$$
\left|\bigcap_{i=1}^tV_i\right|\geq\sum_{i=1}^t|V_i|-(t-1)\left|\bigcup_{i=1}^tV_i\right|.
$$
\end{lem}

Recall the classical stability theorem proved by Erd\H{o}s and Simonovits:
\begin{lem}[\cite{E1966,E1968,Simo1968}]\label{lem3} 
Let \(H\) be a graph with \(\chi \left( H\right)  = r + 1 \geq  3\). For every \(\varepsilon  > 0\), there exist a constant \(\delta > 0\) and an integer \({n}_{0}\) such that if \(G\) is an \(H\)-free graph on \(n \geq  {n}_{0}\) vertices with \(e\left( G\right)  \geq  \left( {1 - \frac{1}{r} - \delta}\right)\cdot \frac{{n}^{2}}{2}\), then \(G\) can be obtained from \({T}_{r}(n)\) by adding and deleting at most \(\varepsilon {n}^{2}\) edges.
\end{lem}

The following result was showed by Amin et al. \cite[Lemma 12]{AFGS}.
\begin{lem}[\cite{AFGS}]\label{lem2.4} 
Let $n,r,q,p$ be integers satisfy $r\geq3,\,n=rq+p\geq r+3,$ and $\,0\leq p\leq r-1.$ Let $G = G_1\vee G_2$ be an $n$-vertex graph such that $G_1\in \mathcal{C}_5^\ast$ and $G_2$ is a complete
$(r-2)$-partite graph. If the graph $G$ has the maximum size, then
$$
G\in\left\{
                                        \begin{array}{ll}
                            \left\{C_5\vee T_{r-2}(n-5)\right\}, & \text{if $q=1,2;$} \\
                       \mathcal{G}_1[n]\cup\mathcal{G}_2[n], &  \text{if $q\geq 3$ and $p=0;$} \\
\mathcal{G}_1[n]\cup\mathcal{G}_2[n]\cup\mathcal{G}_3[n], &  \text{if $q\geq 3$ and $1\leq p\leq r-3;$} \\
\mathcal{G}_2[n]\cup\mathcal{G}_3[n], &  \text{if $q\geq 3$ and $p=r-2;$} \\
\mathcal{G}_3[n], &  \text{if $q\geq 3$ and $p=r-1.$}
                                        \end{array}
                                      \right.
$$
\end{lem}
\section{\normalsize Proof of Theorem~\ref{thm1}}\label{s3}
In this section, we give the proof of Theorem~\ref{thm1}, which determines the maximum number of edges in
non-$r$-partite \({B}_{r,k}\)-free graphs of order $n$, and characterizes all the corresponding extremal graphs. 

Fix integers $r,k,n,p,q$ with $r\geq3,\,k\geq1,\,0\leq p\leq r-1$ and $n=qr+p$ being sufficiently large. Let \(G\) be a graph with maximum size among all non-$r$-partite \({B}_{r,k}\)-free graphs of order \(n\), and denote by $d(v):=d_G(v)$ for $v\in V(G)$. In the following, we are going to progressively refine the structure of $G,$ and show $G$ is isomorphic to one of the graphs presented in Theorem \ref{thm1}(b). Note that if a graph is $K_{r+1}$-free, then it must be \({B}_{r,k}\)-free. By the choice of $G$ and Theorem \ref{thm1.4}(a), one has 
\begin{align}\label{eq:3.1}
e\left( G\right)  \geq  \left( {1 - \frac{1}{r}}\right)  \cdot  \frac{{n}^{2}}{2} - \frac{n}{r} + \frac{p\left( {p + 2}\right) }{2r} - \frac{p}{2} + 1.
\end{align}
 
In the remainder of this section, let \(\varepsilon\) be a fixed constant with \(0<\varepsilon<\frac{1}{36r^8}\).
\begin{lem}\label{lem3.1} 
It holds that
$$
e\left( G\right)  \geq  e\left( {{T}_{r}\left( n\right) }\right)  - \frac{\varepsilon }{2}{n}^{2}.
$$
Further, \(G\) admits a partition \(V\left( G\right)  = {V}_{1} \cup  V_2\cup \cdots  \cup  {V}_{r}\) such that \(\sum_{1 \leq  i < j \leq  r}e\left( {{V}_{i},{V}_{j}}\right)\) attains the maximum, \(\sum_{i =  1}^{r}e\left( {V}_{i}\right)  \leq  \frac{\varepsilon }{2}{n}^{2}\) and for each \(i \in  \left\lbrack  r\right\rbrack\),
$$
\left( {\frac{1}{r} - 2\sqrt{\varepsilon }}\right) n < \left| {V}_{i}\right|  < \left( {\frac{1}{r} + 2\sqrt{\varepsilon }}\right) n.
$$
\end{lem}
\begin{proof}
Note that $n$ is sufficiently large. By \eqref{eq:3.1} and Lemmas \ref{lem2.1}, \ref{lem3}, \(G\) can be obtained from \({T}_{r}\left( n\right)\) by adding or deleting at most \(\frac{\varepsilon }{2}{n}^{2}\) edges. Hence, \(e\left( G\right)  \geq  e\left( {{T}_{r}\left( n\right) }\right)\)  \(- \frac{\varepsilon }{2}{n}^{2}\). Furthermore, there is a partition \(V\left( G\right)  = U_1\cup U_2\cup\cdots \cup U_{r}\) with \(\sum_{i = 1}^{r}e\left( {U}_{i}\right)  \leq  \frac{\varepsilon }{2}{n}^{2}\) and \(\left\lfloor \frac{n}{r}\right\rfloor  \leq  \left| {U}_{i}\right|  \leq  \left\lceil  \frac{n}{r}\right\rceil\) for each \(i \in  \left\lbrack  r\right\rbrack\). Let \(V\left( G\right)  = {V}_{1} \cup V_2\cup \cdots  \cup  {V}_{r}\) be a partition such that \(\sum_{1 \leq  i < j \leq  r}e\left( {{V}_{i},{V}_{j}}\right)\) attains the maximum. Then \(\sum_{i =  1}^{r}e\left( {V}_{i}\right)  \leq  \sum_{i = 1}^{r}e\left( {U}_{i}\right)  \leq  \frac{\varepsilon }{2}{n}^{2}\) and \(\sum_{1 \leq  i < j \leq  r}e\left( {{V}_{i},{V}_{j}}\right)  \geq  \sum_{{1 \leq  i < j \leq  r}}e\left( {{U}_{i},{U}_{j}}\right)  \geq  e\left( {{T}_{r}\left( n\right) }\right)  - \frac{\varepsilon }{2}{n}^{2}\).

Without loss of generality, we assume \(||V_1|-\frac{n}{r}| = \max \{||V_j|-\frac{n}{r}|, j \in [r]\}\). Suppose to the contrary that \(||V_1|-\frac{n}{r}|   \geq  2\sqrt{\varepsilon }n\). Then, 
\begin{align}
e\left( G\right)  &\leq  \sum_{1 \leq  i < j \leq  r}\left| {V}_{i}\right| \left| {V}_{j}\right|  + \sum_{i = 1}^{r}e\left( {V}_{i}\right)\notag\\
&\leq  \left| {V}_{1}\right| \left( {n - \left| {V}_{1}\right| }\right)  + \mathop{\sum }\limits_{{2 \leq  i < j \leq  r}}\left| {V}_{i}\right| \left| {V}_{j}\right|  + \frac{\varepsilon }{2}{n}^{2}\notag\\
&= \left| {V}_{1}\right| \left( {n - \left| {V}_{1}\right| }\right)  + \frac{1}{2}\left[{{\left( \sum_{i = 2}^{r}\left| {V}_{i}\right| \right) }^{2} - \sum_{i = 2}^{r}{\left| {V}_{i}\right| }^{2}}\right]  + \frac{\varepsilon }{2}{n}^{2}\notag\\
&\leq  \left| {V}_{1}\right| \left( {n - \left| {V}_{1}\right| }\right)  + \frac{1}{2}{\left( n - \left| {V}_{1}\right| \right) }^{2} - \frac{1}{2\left( {r - 1}\right) }{\left( n - \left| {V}_{1}\right| \right) }^{2} + \frac{\varepsilon }{2}{n}^{2}\notag\\
%&= \left| {V}_{1}\right| \left( {n - \left| {V}_{1}\right| }\right)  + \frac{r - 2}{2\left( {r - 1}\right) }{\left( n - \left| {V}_{1}\right| \right) }^{2} + \frac{\varepsilon }{2}{n}^{2}\notag\\
%&= \frac{r - 1}{2r}{n}^{2} - \frac{1}{{2r}\left( {r - 1}\right) }{n}^{2} + {\left| {V}_{1} \right\cdot  n - \left|{V}_{1}\right| }^{2} - \frac{r - 2}{r - 1} \cdot  \left| {V}_{1}\right|  \cdot  n + \frac{r - 2}{2\left( {r - 1}\right) }{\left| {V}_{1}\right| }^{2} + \frac{\varepsilon }{2}{n}^{2}\notag\\
&= \frac{r - 1}{2r}{n}^{2} - \frac{1}{{2r}\left( {r - 1}\right) }{n}^{2} + \frac{1}{r - 1}\left| {V}_{1}\right|  \cdot  n - \frac{r}{2\left( {r - 1}\right) }{\left| {V}_{1}\right| }^{2} + \frac{\varepsilon }{2}{n}^{2}\notag\\
%&= \frac{r - 1}{2r}{n}^{2} - \frac{r}{2\left( {r - 1}\right) }\left( {\frac{{n}^{2}}{r} - 2 \cdot  \left| {V}_{1}\right|  \cdot  \frac{n}{r} + {\left| {V}_{1}\right| }^{2}}\right)  + \frac{\varepsilon }{2}{n}^{2}\notag\\
&= \frac{r - 1}{2r}{n}^{2} - \frac{r}{2\left( {r - 1}\right) }{\left( \frac{n}{r} - \left| {V}_{1}\right| \right) }^{2} + \frac{\varepsilon }{2}{n}^{2}\notag\\
&\leq  \frac{r - 1}{2r}{n}^{2} - \frac{r}{2\left( {r - 1}\right) } \cdot  {4\varepsilon }{n}^{2} + \frac{\varepsilon }{2}{n}^{2}\notag\\
&< \frac{r - 1}{2r}{n}^{2} - \frac{3}{2}\varepsilon {n}^{2}. \label{q2}
\end{align}
On the other hand,
$$
e\left( G\right)  \geq  e\left( {{T}_{r}\left( n\right) }\right)  - \frac{\varepsilon }{2}{n}^{2} \geq  \frac{r - 1}{2r}{n}^{2} - \frac{r}{8} - \frac{\varepsilon }{2}{n}^{2} > \frac{r - 1}{2r}{n}^{2} - \varepsilon {n}^{2}
$$
for large enough \(n\), a contradiction to \eqref{q2}.
\end{proof}

\begin{lem}\label{lem3.2} 
Denote
$$
W :=\bigcup_{i = 1}^{r}\left\{{v \in  {V}_{i} :  {d}_{{V}_{i}}\left( v\right)  \geq  3\sqrt{\varepsilon }n}\right\}
$$
and
$$
L :=\left\{  {v \in  V\left( G\right) :  d\left( v\right)  \leq  \left( {1 - \frac{1}{r} - 5\sqrt{\varepsilon}}\right) n}\right\}.
$$
Then \(\left| L\right|  \leq  \sqrt{\varepsilon }n\) and \(W\subseteq L.\)
\end{lem}
\begin{proof}
We first prove the following claims, which establish upper bounds respectively for $|W|$ and $|L|$.
\begin{claim}\label{c1}
\(\left| W\right|  \leq  \frac{1}{3}\sqrt{\varepsilon }n\).
\end{claim}
\begin{proof}[\bf Proof of Claim \ref{c1}] 
It follows from Lemma \ref{lem3.1} that \(\sum_{i = 1}^{r}e\left( {V}_{i}\right)  \leq  \frac{\varepsilon }{2}{n}^{2}\). On the other hand, let \({W}_{i} :={V}_{i} \cap  W\) for all \(i \in  \left\lbrack  r\right\rbrack\). Then
$$
2e\left(V_i\right) =\sum_{u \in  {V}_{i}}{d}_{{V}_{i}}\left( u\right)  \geq \sum_{u \in  {W}_{i}}{d}_{{V}_{i}}\left( u\right)  \geq  \left| {W}_{i}\right|  \cdot  3\sqrt{\varepsilon }n.
$$
Thus,
$$
\frac{\varepsilon}{2}{n}^{2} \geq \sum_{i = 1}^{r}e\left( {V}_{i}\right)  \geq  \left| W\right|  \cdot  \frac{3\sqrt{\varepsilon}}{2}n.
$$
So we obtain $\left| W\right|  \leq  \frac{1}{3}\sqrt{\varepsilon}n$.
\end{proof}
\begin{claim}\label{c2} 
\(\left| L\right|  \leq  \sqrt{\varepsilon }n\).
\end{claim}
\begin{proof}[\bf Proof of Claim~\ref{c2}] Suppose to the contrary that \(\left| L\right|  > \sqrt{\varepsilon }n\). Then there is a subset \(L' \subseteq  L\) with \(\left| L'\right|  = \left\lfloor  {\sqrt{\varepsilon }n}\right\rfloor\). Therefore,
\begin{align*}
e\left( {G\left\lbrack  {V \backslash  {L}'}\right\rbrack}\right)&\geq  e\left( G\right)- \sum_{{v \in  {L}'}}d\left( v\right)\\
&\geq  e\left( {T}_r(n)\right)  - \frac{\varepsilon }{2}{n}^{2} - \sqrt{\varepsilon }n \cdot  \left( {1 - \frac{1}{r} - 5\sqrt{\varepsilon }}\right) n\\
&\geq  \left( {1 - \frac{1}{r}}\right) \frac{{n}^{2}}{2} - \frac{r}{8} - \left( {1 - \frac{1}{r} - \frac{9}{2}\sqrt{\varepsilon }}\right) \sqrt{\varepsilon }{n}^{2}\\
&= \left( {1 - \frac{1}{r}}\right)  \cdot  \frac{{\left( n - \left\lfloor  \sqrt{\varepsilon }n\right\rfloor  \right) }^{2}}{2} + \left( {1 - \frac{1}{r}}\right)  \cdot  n \cdot  \left\lfloor  {\sqrt{\varepsilon }n}\right\rfloor   - \left( {1 - \frac{1}{r}}\right)  \cdot  \frac{{\left\lfloor  \sqrt{\varepsilon }n\right\rfloor  }^{2}}{2}\\
&\,\,\,\,\,\,\,
- \frac{r}{8} - \left( {1 - \frac{1}{r} - \frac{9}{2}\sqrt{\varepsilon }}\right) \sqrt{\varepsilon }{n}^{2}\\
&\geq  \left( {1 - \frac{1}{r}}\right)  \cdot  \frac{{\left(n-\left\lfloor  \sqrt{\varepsilon }n\right\rfloor\right) }^{2}}{2} + \left( {1 - \frac{1}{r}}\right)  \cdot  n \cdot  \left( {\sqrt{\varepsilon }n - 1}\right)  - \left( {1 - \frac{1}{r}}\right)  \cdot  \frac{\varepsilon {n}^{2}}{2}\\
&\,\,\,\,\,\,\,- \frac{r}{8} - \left( {1 - \frac{1}{r}}\right) \sqrt{\varepsilon }{n}^{2} + \frac{9}{2}\varepsilon {n}^{2}\\
&= \left( {1 - \frac{1}{r}}\right)  \cdot  \frac{{\left( n - \left\lfloor  \sqrt{\varepsilon }n\right\rfloor\right) }^{2}}{2} - \left( {1 - \frac{1}{r}}\right) n + \left( {8 + \frac{1}{r}}\right) \frac{\varepsilon {n}^{2}}{2} - \frac{r}{8}\\
&> e\left( {T}_r\left(n - \left\lfloor  \sqrt{\varepsilon }n\right\rfloor\right)\right).
\end{align*}
Note that \({B}_{r,k}\) is a color-critical graph with $\chi({B}_{r,k})=r+1$, and $n$ is sufficiently large. By Theorem \ref{thm1.2}, $e\left( {T}_r\left(n - \left\lfloor  \sqrt{\varepsilon }n\right\rfloor\right)\right)  = \ex\left({n - \left\lfloor  \sqrt{\varepsilon }n\right\rfloor,{B}_{r,k}}\right).$ Hence, \(e\left( {G\left\lbrack  {V \backslash  {L}'}\right\rbrack  }\right)  > \ex\left({n - \left\lfloor  \sqrt{\varepsilon }n\right\rfloor,{B}_{r,k}}\right)\), which implies that \(G\left\lbrack  {V \backslash   L'}\right\rbrack\) and so $G$ contain a copy of \({B}_{r,k}\), a contradiction. 
\end{proof}

By Claim \ref{c2}, the first part of Lemma~\ref{lem3.2} follows immediately.

Next, we prove that \(W \subseteq  L\). Otherwise, there is a vertex \({u}_{0} \in  W \backslash  L\). Without loss of generality, let \({u}_{0} \in  {V}_{1}\). Since \(V\left( G\right)  = {V}_{1} \cup V_2 \cdots  \cup  {V}_{r}\) is the partition such that \(\sum_{{1 \leq  i < j \leq  r}}e\left( {{V}_{i},{V}_{j}}\right)\) attains the maximum, \({d}_{{V}_{1}}\left( {u}_{0}\right)  \leq  {d}_{{V}_{i}}\left( {u}_{0}\right)\) for each \(i \in [r]\backslash\{1\}\). Thus \(d\left( {u}_{0}\right)  \geq  r{d}_{{V}_{1}}\left( {u}_{0}\right)\), that is \({d}_{{V}_{1}}\left( {u}_{0}\right)  \leq  \frac{1}{r}d\left( {u}_{0}\right)\). On the other hand, since \({u}_{0} \notin  L\), we get \(d\left( {u}_{0}\right)  \geq  \left( {1 - \frac{1}{r} - 5\sqrt{\varepsilon }}\right) n\). Thus
\begin{align}
{d}_{{V}_{2}}\left( {u}_{0}\right)  &= d\left( {u}_{0}\right)  - {d}_{{V}_{1}}\left( {u}_{0}\right)  - \sum_{{i = 3}}^{r}{d}_{{V}_{i}}\left( {u}_{0}\right)\notag\\
&\geq  \left( {1 - \frac{1}{r}}\right) d\left( {u}_{0}\right)  - \sum_{{i = 3}}^{r}\left| {V}_{i}\right|\notag\\
&>  \left( {1 - \frac{1}{r}}\right) \left( {1 - \frac{1}{r} - 5\sqrt{\varepsilon }}\right) n - \left( {r - 2}\right) \left( {\frac{1}{r} + 2\sqrt{\varepsilon }}\right)n\notag\\
%&= \left( {1 - \frac{2}{r} + \frac{1}{{r}^{2}} - \frac{r - 1}{r} \cdot  5\sqrt{\varepsilon }}\right) n - \left( {1 - \frac{2}{r} + {2r}\sqrt{\varepsilon } - 4\sqrt{\varepsilon }}\right) n\notag\\
&= \left( {\frac{1}{{r}^{2}} - \left( {{2r} + 1 - \frac{5}{r}}\right) \sqrt{\varepsilon }}\right) n.\label{eq:3.3}
\end{align}
%So we obtain
%\[\varepsilon  < \frac{1}{{16}{r}^{6}}.\]%(3)
Since \(\left| W\right|  \leq  \frac{1}{3}\sqrt{\varepsilon }n,\,\left| L\right|  \leq  \sqrt{\varepsilon }n\). It follows that $d_{V_2\setminus (W\cup L)}(u_0)  > \frac{n}{{r}^{2}} - \left( {{2r} + \frac{7}{3} - \frac{5}{r}}\right) \sqrt{\varepsilon }n$.

On the other hand, note that \(u_0 \in  W,\) and so \(d_{V_1}(u_0)  \geq  3\sqrt{\varepsilon}n>\left|L\cup W\right|\). Therefore, we can choose a vertex \(u_1 \in  N_{V_1}(u_0)  \backslash (W \cup  L)\). Then
\begin{align}
{d}_{{V}_{2}}\left( {u}_{1}\right)  &= d\left( {u}_{1}\right)  - {d}_{{V}_{1}}\left( {u}_{1}\right)  - \sum_{{i = 3}}^{r}{d}_{{V}_{i}}\left( {u}_{1}\right)\notag\\
&> \left( {1 - \frac{1}{r} - 5\sqrt{\varepsilon }}\right) n - 3\sqrt{\varepsilon }n - \sum_{{i = 3}}^{r}\left|{V}_{i}\right|\notag\\
&> \left( {1 - \frac{1}{r} - 8\sqrt{\varepsilon }}\right) n - \left( {r - 2}\right) \left( {\frac{1}{r} + 2\sqrt{\varepsilon }}\right) n\notag\\
&= \frac{n}{r} - \left( {{2r} + 4}\right) \sqrt{\varepsilon }n. \label{eq:3.4}
\end{align}
Therefore, by Lemmas \ref{lem2.2}, \ref{lem3.1} and \eqref{eq:3.3}-\eqref{eq:3.4},
\begin{align*}
\left| {{N}_{{V}_{2}}\left( {u}_{0}\right)  \cap  {N}_{{V}_{2}}\left( {u}_{1}\right) }\right|  &\geq  \left| {{N}_{{V}_{2}}\left( {u}_{0}\right) }\right|  + \left| {{N}_{{V}_{2}}\left( {u}_{1}\right) }\right|  - \left| {{N}_{{V}_{2}}\left( {u}_{0}\right)  \cup  {N}_{{V}_{2}}\left( {u}_{1}\right) }\right|\\
&> \frac{n}{{r}^{2}} - \left( {{2r} + 1 - \frac{5}{r}}\right) \sqrt{\varepsilon }n + \frac{n}{r} - \left( {{2r} + 4}\right) \sqrt{\varepsilon }n - \left( {\frac{1}{r} + 2\sqrt{\varepsilon }}\right) n\\
&> \frac{n}{{r}^{2}} - \left( {{4r} + 7}\right) \sqrt{\varepsilon}n.
\end{align*}
Since \(\left| W\right|  \leq  \frac{1}{3}\sqrt{\varepsilon }n,\,\left| L\right|  \leq  \sqrt{\varepsilon }n\), one has
\begin{align*}
\left| \left( {{N}_{{V}_{2}}\left( {u}_{0}\right)  \cap  {N}_{{V}_{2}}\left( {u}_{1}\right) }\right) \setminus \left({W \cup  L}\right) \right| &>\frac{n}{{r}^{2}} - \left( {{4r} + 7}\right) \sqrt{\varepsilon }n - \frac{4}{3}\sqrt{\varepsilon }n\\
&=\frac{n}{{r}^{2}} - \left( {{4r} + \frac{25}{3}}\right) \sqrt{\varepsilon }n > 0.
\end{align*}
%Hence,
%\[\varepsilon  < \frac{1}{{64}{r}^{6}}.\]

Hence, there is a vertex \({u}_{2}\) in \({V}_{2} \setminus  \left( {W \cup  L}\right)\) being adjacent to both \({u}_{0}\) and \({u}_{1}\). For an integer \(s\) with \(2 \leq  s \leq  r - 1\), assume that for any \(1 \leq  i \leq  s\), there is a vertex \(u_i\in V_i \setminus  \left( {W\cup L}\right)\) such that \(\left\{  {{u}_{0},{u}_{1},\ldots,{u}_{s}}\right\}\) is a clique. We next consider the number of common neighbors of these vertices in \({V}_{s + 1} \backslash  \left( W\cup L\right)\). Similarly, by \eqref{eq:3.3} and \eqref{eq:3.4}, we get that
$$
{d}_{{V}_{s + 1}}\left( {u}_{0}\right) > \frac{n}{{r}^{2}} - \left( {{2r} + 1 - \frac{5}{r}}\right) \sqrt{\varepsilon}n,
$$
and for each $i\in [s]$,
$$
{d}_{{V}_{s + 1}}\left( {u}_{i}\right)  > \frac{n}{r} - \left( {{2r} + 4}\right) \sqrt{\varepsilon }n.
$$
Together with Lemmas \ref{lem2.2} and \ref{lem3.1}, one gets
\begin{align*}
\left|N_{V_{s + 1}}(u_0)\bigcap \left(\bigcap_{i = 1}^sN_{V_{s + 1}}(u_i)\right)\right|  &\geq  {d}_{{V}_{s + 1}}\left( {u}_{0}\right)  + \mathop{\sum }\limits_{{i = 1}}^{s}{d}_{{V}_{s + 1}}\left( {u}_{i}\right)  - s \cdot  \left| {V}_{s + 1}\right|\\
&>\frac{n}{{r}^{2}} - \left( {{2r} + 1 - \frac{5}{r}}\right) \sqrt{\varepsilon}n + s \cdot  \left( {\frac{n}{r} - \left( {{2r} + 4}\right) \sqrt{\varepsilon }n}\right)  \\
&\,\,\,\,\,\,\,-s \cdot  \left( {\frac{1}{r} + 2\sqrt{\varepsilon }}\right) n\\
&>\frac{n}{{r}^{2}} - \left( {{2sr} + {2r+1} + {6s}}\right) \sqrt{\varepsilon }n. 
\end{align*}
%So, 
%$$
 %\varepsilon  < \frac{1}{{36}{r}^{8}}.
%$$
Hence 
$$\left|N_{V_{s + 1}}(u_0)  \bigcap  \left(\bigcap_{i = 1}^{s}N_{V_{s + 1}}(u_i)\right)\bigg\backslash \left({W \cup  L}\right) \right| >  \frac{n}{{r}^{2}} - \left( {{2sr} + {2r} + {6s} + \frac{7}{3}}\right) \sqrt{\varepsilon }n > k.
$$ 
That is to say, \({u}_{0},{u}_{1},\ldots ,{u}_{s}\) have at least \(k\) common neighbors in \({V}_{s + 1} \backslash  \left( {W \cup  L}\right)\). Thus, for each \(2\leq i\leq r-1\), there exists a vertex \(u_i \in  V_i\setminus (W \cup  L)\) such that \(\left\{  {{u}_{0},{u}_{1},\ldots ,{u}_{r - 1}}\right\}\) is a clique, and they have \(k\) common neighbors \({u}_{r,1},{u}_{r,2},\ldots ,{u}_{r,k}\) in \({V}_{r}\). Now \(G\left\lbrack  \left\{  {{u}_{0},{u}_{1},\ldots ,{u}_{r - 1},{u}_{r,1},{u}_{r,2},\ldots ,{u}_{r,k}}\right\}  \right\rbrack\) is isomorphic to \({B}_{r,k}\), a contradiction. Hence, $W\subseteq L$.
\end{proof} 
\begin{lem}\label{lem3.3} 
\(\chi(G)  \geq  r + 1\).
\end{lem} 
\begin{proof}
Denote by \(s = \chi \left( G\right)\). Since \(G\) is non-$r$-partite with \(s \neq  r\). If \(s \leq  r - 1\), then by Lemma \ref{lem2.1}, \(e\left( G\right)  \leq  e\left( {{T}_{s}\left( n\right) }\right)  \leq  \left( {1 - \frac{1}{s}}\right) \frac{{n}^{2}}{2} < \left( {1 - \frac{1}{r}}\right)  \cdot  \frac{{n}^{2}}{2} - \frac{n}{r} + \frac{p\left( {p + 2}\right) }{2r} - \frac{p}{2} + 1\), a contradiction to \eqref{eq:3.1}. Therefore, \(\chi \left( G\right)  \geq  r + 1\). 
\end{proof}

\begin{lem}\label{lem3.4} 
For each \(i \in  \left\lbrack  r\right\rbrack  ,e\left( {{V}_{i} \backslash  L}\right)  = 0\).
\end{lem}
\begin{proof}
Suppose to the contrary that there is an \({i}_{0} \in  \left\lbrack  r\right\rbrack\) such that
\(e\left(V_{i_0}\setminus L\right)  \geq  1\). Without loss of generality, we may assume that \(e\left(V_1\setminus L\right)  \geq  1\). Let
\({u}_{0}{u}_{1}\) be an edge in \(G\left\lbrack V_1\setminus  L\right\rbrack\). Now \({u}_{0},{u}_{1} \notin  L\), and so by Lemma \ref{lem3.2}, \({u}_{0},{u}_{1} \notin  W\).
Hence, $d\left( {u}_{0}\right)  > \left( {1 - \frac{1}{r} - 5\sqrt{\varepsilon}}\right) n$ and ${d}_{{V}_{1}}\left({u}_{0}\right)  < 3\sqrt{\varepsilon }n.$ Together with Lemma \ref{lem3.1}, one has
\begin{align}
{d}_{{V}_{2}}\left( {u}_{0}\right)  &= d\left( {u}_{0}\right)  - {d}_{{V}_{1}}\left( {u}_{0}\right)  - \mathop{\sum }\limits_{{i = 3}}^{r}{d}_{{V}_{i}}\left( {u}_{0}\right)\notag\\
&> \left( {1 - \frac{1}{r} - 5\sqrt{\varepsilon }}\right) n - 3\sqrt{\varepsilon }n - \left( {r - 2}\right) \left( {\frac{1}{r} + 2\sqrt{\varepsilon }}\right) n\notag\\
&= \left( {\frac{1}{r} - \left( {{2r} + 4}\right) \sqrt{\varepsilon }}\right) n.\label{q2.7}%\tag{5}
\end{align}
Similarly, \({d}_{{V}_{2}}\left( {u}_{1}\right)  \geq  \left( {\frac{1}{r} - \left( {{2r} + 4}\right) \sqrt{\varepsilon }}\right) n.\) Together with Lemmas \ref{lem2.2} and \ref{lem3.1}, one has
\begin{align*}
\left| {{N}_{{V}_{2}}\left( {u}_{0}\right)  \cap  {N}_{{V}_{2}}\left( {u}_{1}\right) }\right|  &\geq  {d}_{{V}_{2}}\left( {u}_{0}\right)  + {d}_{{V}_{2}}\left( {u}_{1}\right)  - \left| {{N}_{{V}_{2}}\left( {u}_{0}\right)  \cup  {N}_{{V}_{2}}\left( {u}_{1}\right) }\right|\\
&> \left( {\frac{2}{r} - 2\left( {{2r} + 4}\right) \sqrt{\varepsilon }}\right) n - \left( {\frac{1}{r} + 2\sqrt{\varepsilon }}\right) n\\
&= \left( {\frac{1}{r} - \left( {{4r} + {10}}\right) \sqrt{\varepsilon}}\right)n.
\end{align*}
%Hence,
%\[
%\varepsilon  < \frac{1}{{49}{r}^{4}},
%\]
Note that $|L|\leq\sqrt{\varepsilon }n$ (see Lemma \ref{lem3.2}). Hence, 
$$\left| {\left( {{N}_{{V}_2}\left( {u}_{0}\right)  \cap  {N}_{{V}_2}\left( {u}_{1}\right) }\right)  \setminus  L}\right|  >  \left( {\frac{1}{r} - \left( {{4r} + {11}}\right) \sqrt{\varepsilon }}\right) n > 0.
$$ 
That is to say, there is a vertex
\({u}_{2}\) in \(V_2\setminus L\) being adjacent to both \({u}_{0}\) and \({u}_{1}\). For an integer \(s\) with \(2 \leq  s \leq  r-1\), 
assume that for any \(1 \leq  i \leq  s\), there is a vertex \({u}_{i} \in  {V}_{i} \setminus  L\) such that \(\{{{u}_{0},{u}_{1}},\ldots, u_s\}\) is a clique. 

We next consider the number of common neighbors of these vertices in \({V}_{s + 1} \setminus  L\). Similarly, by \eqref{q2.7} we get 
$$
{d}_{{V}_{s + 1}}\left( {u}_{i}\right)  > \left( {\frac{1}{r} - \left( {{2r} + 4}\right) \sqrt{\varepsilon }}\right) n
$$
for each $i\in \{ 0,1,\ldots, s\}$. Hence, by Lemmas \ref{lem2.2} and \ref{lem3.1},
\begin{align*}
\left|\bigcap_{{i = 0}}^{s}{N}_{{V}_{s + 1}}\left( {u}_{i}\right)\right|  &\geq  \sum_{{i = 0}}^{s}d_{V_{s + 1}}\left( {u}_{i}\right)  - s \cdot  \left| {V}_{s + 1}\right|
\\
&> \left( {s + 1}\right) \left( {\frac{1}{r} - \left( {{2r} + 4}\right) \sqrt{\varepsilon }}\right) n - s\left( {\frac{1}{r} + 2\sqrt{\varepsilon }}\right) n\\
&= \frac{n}{r} - \left( {{2rs} + {6s} + {2r} + 4}\right) \sqrt{\varepsilon }n.
\end{align*}
%$\Rightarrow \varepsilon  < \frac{1}{{36}{r}^{6}}$. 
Therefore, \({u}_{0},{u}_{1},\ldots ,{u}_{s}\) have at least \(\left|\bigcap_{{i = 0}}^{s}{N}_{{V}_{s + 1}}\left( {u}_{i}\right) \right|  - \left| L\right|  > \frac{n}{r} - ({2rs} + {6s} + {2r}\) +5) \(\sqrt{\varepsilon }n \geq  k\) common neighbors in \({V}_{s + 1}\setminus L\). Thus, for each \(2\leq i\leq r - 1\), there exists \({u}_{i} \in  {V}_{i} \setminus L\) such that \(\left\{  {{u}_{0},{u}_{1},\ldots ,{u}_{r - 1}}\right\}\) is a clique, and they have \(k\) common neighbors \({u}_{r,1},{u}_{r,2},\ldots ,{u}_{r,k}\) in \({V}_{r}\). Now \(G\left\lbrack  \left\{{{u}_{0},{u}_{1},\ldots ,{u}_{r-1},{u}_{r,1},{u}_{r,2},\ldots ,{u}_{r,k}}\right\}  \right\rbrack\) is isomorphic to \({B}_{r,k}\), a contradiction. 
\end{proof}

Note that by Lemma \ref{lem3.3}, \(\chi(G)  \geq  r + 1\). It follows from Lemma \ref{lem3.4} that $L\not=\emptyset$. Our next lemma shows $\left| L\right|  \leq  2.$

\begin{lem}\label{lem3.5} 
\(1 \leq  \left| L\right|  \leq  2\).
\end{lem}
\begin{proof}Since \(L  \neq  \emptyset\), one has \(\left| L\right|  \geq  1\). Suppose to the contrary that \(\left| L\right|  \geq  3\). Since \(\chi(G)  \geq  r + 1\), \(\sum_{{i = 1}}^{r}e\left(V_{i}\right)  \neq  0\). Without loss of generality, we may assume \(e\left(V_{1}\right)  \neq  0\). Let \({uv} \in  E\left(V_{1}\right)\), by Lemma~\ref{lem3.4}, at least one of \(u\) and \(v\) is in $L$. Since \(\left| L\right| \geq 3\), one may choose a vertex \(w \in  L \backslash  \{ u,v\}\) and construct a new graph \(G' = G - \left\{  {w{w}^{\prime } : {w}^{\prime } \in  {N}_{G}\left( w\right) }\right\}   + \left\{  {w{w}^{\prime \prime } : {w}^{\prime \prime } \in  V\left( G\right)  \backslash  \left(V_1\cup L\right) }\right\}\). Then \({G}^{\prime }\) is \({B}_{r,k}\)-free and non-$r$-partite. However,
\begin{align*}
e({G}^{\prime })&= e\left( G\right)  - {d}_{G}\left( w\right)  + \left|{V\left( G\right)  \backslash  \left(V_1\cup L\right) }\right|\\
&\geq e(G)  - \left( {1 - \frac{1}{r} - 5\sqrt{\varepsilon }}\right) n + n - \left( {\frac{1}{r} + 2\sqrt{\varepsilon }}\right) n - \sqrt{\varepsilon }n\\
&= e\left( G\right)  + 2\sqrt{\varepsilon }n\\
&> e\left( G\right),
\end{align*}
a contradiction to the choice of \(G\). Hence, \(\left| L\right|  \leq  2\).
\end{proof}
\begin{lem}\label{lem3.6}
\(\sum_{{i = 1}}^{r}e\left( {V}_{i}\right) = 1\).
\end{lem}
\begin{proof}
Since \(\chi(G)  \geq  r + 1\) (see Lemma \ref{lem3.3}), one has \(\sum_{{i = 1}}^{r}e\left( {V}_{i}\right)  \geq  1\). Suppose to the contrary that \(\sum_{{i = 1}}^{r}e\left(V_i\right)  \geq  2\). Without loss of generality, we may assume 
\(e\left(V_{1}\right)  \geq  e\left(V_{2}\right)  \geq  e\left(V_{3}\right)  \geq  \cdots  \geq  e\left(V_{r}\right)\). Then by Lemmas \ref{lem3.4} and \ref{lem3.5}, \(e\left(V_{3}\right)  = e\left(V_{4}\right)  = \cdots  = e\left(V_{r}\right)  = 0\), and so \(e\left( {V}_{1}\right)  + e\left( {V}_{2}\right)  = \sum_{{i = 1}}^{r}e\left( {V}_{i}\right)  \geq  2\). We complete the proof of this lemma by
considering the following two possible cases.

{\bf Case 1.}\ \(e\left( {V}_{2}\right)  = 0\). In this case, \(e\left( {V}_{1}\right)  \geq  2\). We claim \(L \subseteq  {V}_{1}\). Suppose to the contrary that there is a vertex \(u \in  L \backslash  {V}_{1}\). Then construct \({G}^{\prime } = G - \{ {uw} :  w \in  {N}_{G}\left( u\right) \}  + \left\{ {uw'} :  w' \in  \left(\bigcup_{{i =  2}}^r{V}_{i}\right)  \setminus  \{u\}\right\}\). Since \(e\left( {V}_{1}\right)  \geq  2,\,\sum_{{i = 2}}^{r}e\left( {V}_{i}\right)=0,\) the graph \({G}^{\prime }\) is \({B}_{r,k}\)-free non-$r$-partite. However,
\begin{align*}
e\left( {G}^{\prime }\right)  - e\left( G\right)  &= n - \left| {V}_{1}\right|  - 1 - {d}_{G}\left( u\right)\\
&> n - \left( {\frac{1}{r} + 2\sqrt{\varepsilon }}\right) n - 1 - \left( {1 - \frac{1}{r} - 5\sqrt{\varepsilon }}\right) n\\
&= 3\sqrt{\varepsilon }n - 1 > 0,
\end{align*}
a contradiction to the choice of \(G\). Hence \(L \subseteq  {V}_{1}\).

By Lemma \ref{lem3.5}, \(1 \leq  \left| L\right|  \leq  2\). If \(\left| L\right|  = 1\), let \(V_1\cap L = \{u\}\). Then there are \(e\left(V_1\right)  \geq  2\) vertices in \({V}_{1} \setminus  L\) adjacent to \(u\). 
\begin{claim}\label{claim3}
For each $s\in [r]$ and each vertex in $N_{V_s}(u)$, there is at most one vertex in \(\bigcup_{{i = 1,i \neq  s}}^{r}{V}_{i}\) not being adjacent to it.
\end{claim}
\begin{proof}[\bf Proof of Claim~\ref{claim3}] 
Suppose to the contrary that there is an $s\in [r]$ and a vertex $u_s\in N_{V_s}(u)$ such that there are at least two vertices in \(\bigcup_{{i = 1,i \neq  s}}^{r}{V}_{i}\) not being adjacent to $u_s.$ Then construct \({G}^{\prime } = G - u{u}_s + \{{{u}_sv :  v \in  (\bigcup_{{i = 1,i \neq  s}}^{r}{V}_{i})  \backslash  {N}_{G}\left( {u}_s\right) }\}\). Clearly, $G'$ is non-$r$-partite $B_{r,k}$-free. However, 
$$
e\left( {G}^{\prime }\right)  - e\left( G\right)  = \left| \left( {\bigcup_{{i = 1,i \neq  s}}^{r}{V}_{i}}\right)\bigg\backslash {{N}_{G}\left( {u}_s\right) }\right|  - 1 > 0,
$$ 
a contradiction to the choice of $G$. 
\end{proof}

Recall that $n=qr+p,$ where $0\leq p\leq r-1.$
\begin{claim}\label{claim4}
\(d\left( u\right)  \geq  \left( {1 - \frac{2}{r}}\right) n + \frac{(p+1)^2}{2r} - \frac{p-1}{2}\).
\end{claim}
\begin{proof}[\bf Proof of Claim~\ref{claim4}]
Note that by \eqref{eq:3.1}, $e\left( G\right)  \geq  \left( {1 - \frac{1}{r}}\right)  \cdot  \frac{{n}^{2}}{2} - \frac{n}{r} + \frac{p\left( {p + 2}\right) }{2r} - \frac{p}{2} + 1.$ Since $G- \{u\}$ is an \(r\)-partite graph, by Lemma \ref{lem2.1}, $e\left({G-u}\right)  \leq  e\left(T_r(n - 1)\right)  \leq  \left(1 - \frac{1}{r}\right)  \cdot  \frac{( n - 1)^2}{2}.$ Thus
\begin{equation}\label{eq:3.6}
\left( {1 - \frac{1}{r}}\right)  \cdot  \frac{{n}^{2}}{2} - \frac{n}{r} + \frac{p\left( {p + 2}\right) }{2r} - \frac{p}{2} + 1 \leq  e\left( G\right)  = d\left( u\right)  + e\left( {G - u}\right)  \leq  d\left( u\right)  + \left(1 - \frac{1}{r}\right)  \cdot  \frac{{n}^{2} - {2n} + 1}{2},%\left( 6\right)
\end{equation}
which gives \(d\left( u\right)  \geq  \left( {1 - \frac{2}{r}}\right) n + \frac{(p+1)^2}{2r} - \frac{p-1}{2}\).
\end{proof}
Recall that \(V(G)  = V_1 \cup  \cdots  \cup  V_r\) is a partition maximizing \(\sum_{{1 \leq  i < j \leq  r}}e(V_i,V_j)\). It holds that \(d_{V_i}(u)  \geq  d_{V_1}(u)  = e(V_1)\) for each \(i \in  \{2,\ldots, r\}\). Without loss of generality, we may assume \(d_{V_2}(u)  \leq 
d_{V_3}(u)  \leq  \cdots  \leq  d_{V_r}(u)\).
\begin{claim}\label{claim5}
\({d}_{{V}_{3}}\left( u\right)  \geq  \frac{n}{4r}\), and so \({d}_{{V}_{s}}\left( u\right)  \geq  \frac{n}{4r}\) for all $3\leq s\leq r.$
\end{claim}
\begin{proof}[\bf Proof of Claim \ref{claim5}] Suppose to the contrary that \({d}_{{V}_{3}}\left( u\right)  < \frac{n}{4r}\), then 
\(d_{V_1}(u)  \leq  d_{V_2}(u)  < \frac{n}{4r}\). It follows that
\begin{align*}
d\left( u\right)  = \sum_{{i = 1}}^{r}d_{V_i}\left( u\right)  &= d_{V_1}\left( u\right)  + d_{V_2}\left( u\right)  + d_{V_3}\left( u\right)  + \sum_{{i = 4}}^{r}d_{V_i}\left( u\right)\\
&< \frac{3n}{4r} + \sum_{{i = 4}}^{r}\left| {V}_{i}\right|\\
&= \frac{3n}{4r} + n -\sum_{{i = 1}}^{3}\left| {V}_{i}\right|\\
&<  \frac{3n}{4r} + n - 3\left( {\frac{1}{r} - 2\sqrt{\varepsilon }}\right) n\\
%&= \frac{3n}{4r} + \left( {1 - \frac{3}{r} + 6\sqrt{\varepsilon }}\right) n\\
&= \left( {1 - \frac{9}{4r} + 6\sqrt{\varepsilon}}\right) n\\
&< \left( {1 - \frac{2}{r}}\right) n + \frac{(p+1)^2}{2r} - \frac{p-1}{2},
\end{align*}
a contradiction to Claim \ref{claim4}.
\end{proof}

Take a vertex \(u_1\in N_{V_1}(u)\), then by Claim \ref{claim3}, \(d_{V_i}(u_1)  \geq  \left| {V}_{i}\right|  - 1\) for each \(i \in  \{ 2,\ldots ,r\}\). Since \(d_{V_2}( u)  \geq  d_{V_1}( u)  = e(V_1)  \geq  2\), there is at least one vertex in \({V}_{2}\) being adjacent to both \(u\) and \({u}_{1}\). For an integer \(s\) with \(2 \leq  s \leq  r - 1\), assume that for any \(1 \leq  i \leq  s\), there is a vertex \(u_i\in V_i\) such that \(\{u, u_1,\ldots, u_s\}\) is a clique.
We consider the number of common neighbors of these vertices in \(V_{s+1}\). By Lemma \ref{lem2.2}, Claims \ref{claim3} and \ref{claim5},
\begin{align*}
\left|N_{V_{s + 1}}(u) \bigcap \left( {\bigcap_{{i = 1}}^{s}{N}_{{V}_{s + 1}}\left( {u}_{i}\right) }\right)\right|  &\geq  \left| {{N}_{{V}_{s + 1}}\left( u\right) }\right|  + \mathop{\sum }\limits_{{i = 1}}^{s}\left| {{N}_{{V}_{s + 1}}\left( {u}_{i}\right) }\right|  - s \cdot  \left|V_{s + 1}\right|\\
&\geq  \frac{n}{4r} + s \cdot  \left( {\left| {V}_{s + 1}\right|  - 1}\right)  - s \cdot  \left| {V}_{s + 1}\right|\\
&= \frac{n}{4r} - s \\
&\geq  k.
\end{align*}
That is to say, $u, u_1, \ldots, u_s$ have at least $k$ common neighbors in $V_{s+1}$. Thus, 
for each \(2\leq i\leq r - 1\), there exists \({u}_{i} \in  {V}_{i}\) such that \(\left\{  {u,{u}_{1},\ldots ,{u}_{r - 1}}\right\}\) is a clique, and they have \(k\) common neighbors \({u}_{r,1},{u}_{r,2},\ldots ,{u}_{r,k}\) in \({V}_{r}\). Now
\(G\left\lbrack  {\{ u,{u}_{1},\ldots ,{u}_{r - 1},{u}_{r,1},{u}_{r,2},\ldots ,{u}_{r,k}\} }\right\rbrack\) is isomorphic to \({B}_{r,k}\), a contradiction.

If \(|L| = 2\), then since \(e(V_1)  \geq  2\), there is a vertex \(u\) in \(L\) such that \(e(V_1 \setminus \{u\})\geq 1\). Construct \(G'= G - \left\{  {{uw}:w \in  {N}_{G}\left( u\right) }\right\}   + \left\{  {{uw'}:w' \in \bigcup_{{i = 2}}^{r}{V}_{i}}\right\}\). Then \({G}^{\prime }\) is non-$r$-partite \({B}_{r,k}\)-free. However, 
\begin{align*}
e\left( {G}^{\prime }\right)  - e\left( G\right)  &= \left( {n - \left| {V}_{1}\right| }\right)  - {d}_{G}\left( u\right)\\
&> n - \left( {\frac{1}{r} + 2\sqrt{\varepsilon }}\right) n - \left( {1 - \frac{1}{r} - 5\sqrt{\varepsilon }}\right) n\\
&= 3\sqrt{\varepsilon }n > 0,
\end{align*}
a contradiction to the choice of \(G\).

{\bf Case 2.}\ \(e(V_2)  \geq  1\). Then \(e(V_1)  \geq  e(V_2)  \geq  1\). By Lemmas \ref{lem3.4} and \ref{lem3.5}, \(|V_1 \cap  L|  = |V_2 \cap  L|  = 1\).

Let \(V_2 \cap  L = \{u\}\), construct \(G' = G - \left\{{uw} :  w \in  {N}_{G}\left( u\right)\right\}  + \left\{{uw'} :  w' \in  \left( {{V}_{1} \setminus  L}\right)  \bigcup  \left({\bigcup_{{i = 3}}^{r}{V}_{i}}\right)\right\}\). Then \({G}^{\prime }\) is non-$r$-partite \({B}_{r,k}\)-free. However,
\begin{align*}
e(G')  - e(G)  &= \left( {n - \left| {V}_{2}\right|  - 1}\right)  - {d}_{G}\left( u\right)\\
&> n - \left( {\frac{1}{r} + 2\sqrt{\varepsilon }}\right) n - 1 - \left( {1 - \frac{1}{r} - 5\sqrt{\varepsilon }}\right) n\\
&= 3\sqrt{\varepsilon }n - 1 > 0,
\end{align*}
a contradiction to the choice of \(G\). Therefore, \(\sum_{{i = 1}}^{r}e\left( {V}_{i}\right)  = 1\).
\end{proof}

By Lemma \ref{lem3.6}, we may assume \(e(V_1)= 1\) and \(e(V_i)= 0\) for each \(i \in  \{2,\ldots, r\}\). Let \(uv\) be the unique edge in \(G\left\lbrack  {V}_{1}\right\rbrack\).
\begin{lem}\label{lem3.7}
There is an $i\in \{2, 3, \ldots, r\}$ such that \({N}_{V_i}(u)  \cap  {N}_{V_i}(v)  = \emptyset\).
\end{lem}
\begin{proof}
Suppose to the contrary that for each $i\in \{ 2,3,\ldots ,r\}$, \({N}_{V_i}(u)  \cap  {N}_{V_i}(v)  \not= \emptyset\). Note that both $G-u$ and $G-v$ are $r$-partite graphs. By \eqref{eq:3.6}, \(\min \{d(u), d(v)\}  \geq  \left( {1 - \frac{2}{r}}\right) n + \frac{(p+1)^2}{2r} - \frac{p-1}{2}\). Similarly, notice also that \(G - \{ u,v\}\) is an $r$-partite graph. By Lemma \ref{lem2.1}, \(e(G - \{u,v\})  \leq  e ({T}_{r}(n - 2))  \leq  (1 - \frac{1}{r})  \cdot  \frac{(n - 2)^2}{2}.\) Together with \eqref{eq:3.1}, one may see
\begin{align*}
\left( {1 - \frac{1}{r}}\right)  \cdot  \frac{{n}^{2}}{2} - \frac{n}{r} + \frac{p\left( {p + 2}\right) }{2r} - \frac{p}{2} + 1 &\leq  e \left( G\right)  = e\left( {G - \{ u,v\} }\right)  + d\left( u\right)  + d\left( {v }\right)- 1\\
&\leq  \left( {1 - \frac{1}{r}}\right)  \cdot  \frac{{\left( n - 2\right) }^{2}}{2} + d\left( u\right)  + d\left( {v}\right)  - 1,
\end{align*}
which gives \(d\left( u\right)  + d\left( v\right)  \geq  \left( {2 - \frac{3}{r}}\right) n + \frac{p^2+2p+4}{2r} - \frac{p}{2}\). 

Now
\begin{align}
\left| {N\left( u\right)  \cap  N\left( v\right) }\right|  &= d\left( u\right)  + d\left( v\right)  - \left| {N\left( u\right)  \cup  N\left( v\right) }\right|\notag\\
&\geq  d\left( u\right)  + d\left( v\right)  - 2 - \sum_{{i = 2}}^{r}\left| {V}_{i}\right|\notag\\
&= d\left( u\right)  + d\left( v\right)  - 2 - \left( {n - \left| {V}_{1}\right| }\right)\notag\\
&> \left( {2 - \frac{3}{r}}\right) n + \frac{p^2+2p+4}{2r} - \frac{p}{2} - 2 - n + \left( {\frac{1}{r} - 2\sqrt{\varepsilon }}\right) n\notag\\
&= \left( {1 - \frac{2}{r} - 2\sqrt{\varepsilon }}\right) n + \frac{p^2+2p+4}{2r} - \frac{p}{2} - 2.\label{eq:3.7}
\end{align}
Note that \({N}_{{V}_{1}}\left( u\right)  \cap  {N}_{{V}_{1}}\left( v\right)  = \emptyset\). Without loss of generality, we may assume
$$\left| {{N}_{{V}_{2}}\left( u\right)  \cap  {N}_{{V}_{2}}\left( v\right) }\right|  \leq  \left| {{N}_{{V}_3}\left( u\right)  \cap  {N}_{{V}_{3}}\left( v\right) }\right|  \leq  \cdots  \leq  \left| {{N}_{{V}_{r}}\left( u\right)  \cap  {N}_{{V}_{r}}\left( v\right) }\right|.
$$
Then \(\left| {{N}_{{V}_{2}}\left( u\right)  \cap  {N}_{{V}_{2}}\left( v\right) }\right|  \geq  1\).
\begin{claim}\label{Claim 4}
\(\left| {{N}_{{V}_{3}}\left( u\right)  \cap  {N}_{{V}_{3}}\left( v\right) }\right|  \geq  \frac{n}{4r}\), and so \(\left| {{N}_{{V}_{s}}\left( u\right)  \cap  {N}_{{V}_{s}}\left( v\right) }\right|  \geq  \frac{n}{4r}\) for all $3\leq s\leq r.$
\end{claim}
\begin{proof}[\bf Proof of Claim \ref{Claim 4}] Suppose to the contrary that \(\left| {{N}_{{V}_{3}}\left( u\right)  \cap  {N}_{{V}_{3}}\left( v\right) }\right|  < \frac{n}{4r}\). Then $$
\left| {{N}_{{V}_{2}}\left( u\right)  \cap  {N}_{{V}_{2}}\left( v\right) }\right|  \leq  \left| {{N}_{{V}_3}\left( u\right)  \cap  {N}_{{V}_3}\left( v\right) }\right|  < \frac{n}{4r}.
$$ 
Hence
\begin{align*}
\left| {N\left( u\right)  \cap  N\left( v\right) }\right|  &= \sum_{{i = 2}}^{r}\left| {{N}_{{V}_{i}}\left( u\right)  \cap  {N}_{{V}_{i}}\left( v\right) }\right|\\
&= \left| {{N}_{{V}_{2}}\left( u\right) \cap{N}_{{V}_{2}}\left( v\right) }\right|  + \left| {{N}_{{V}_{3}}\left( u\right) \cap{N}_{{V}_{3}}\left( v\right) }\right|  + \sum_{{i = 4}}^{r}\left| {{N}_{{V}_{i}}\left( u\right) \cap{N}_{{V}_{i}}\left( v\right) }\right|\\
&< \frac{n}{2r} + \sum_{{i = 4}}^{r}\left| {V}_{i}\right|\\
&= \frac{n}{2r} + \left( {n - \left| {V}_{1}\right|  - \left| {V}_{2}\right|  - \left| {V}_{3}\right| }\right)\\
&< \frac{n}{2r} + n - 3 \cdot  \left( {\frac{1}{r} - 2\sqrt{\varepsilon }}\right) n\\
&= \left( {1 - \frac{5}{2r} + 6\sqrt{\varepsilon }}\right) n\\
&< \left( {1 - \frac{2}{r} - 2\sqrt{\varepsilon }}\right) n + \frac{p^2+2p+4}{2r} - \frac{p}{2} - 2,
\end{align*}
a contradiction to \eqref{eq:3.7}. Therefore, \(\left| {{N}_{{V}_{3}}\left( u\right)  \cap  {N}_{{V}_{3}}\left( v\right) }\right|  \geq  \frac{n}{4r}\).
\end{proof}

\begin{claim}\label{Claim 5}
For each \(i \in  \{2,3,\ldots ,r\}\) and each \(w \in  {V}_{i}\), if \(w \in  N(u)  \cap  N(v)\), then there is at most one vertex in \(V(G)  \setminus  V_i\) not adjacent to \(w\).
\end{claim}
\begin{proof}[\bf Proof of Claim \ref{Claim 5}] Suppose to the contrary that there is some $i\in \{2,3,\ldots,r\}$ and some vertex $w\in  {N}_{{V}_{i}}\left( u\right)  \cap  {N}_{{V}_{i}}\left( v\right)$ such that \(w\) is not adjacent to at least two vertices in \(V\left( G\right)  \backslash  {V}_{i}\). Then construct \({G}^{\prime } = G - \left\{  {w{w}^{\prime } :  {w}^{\prime } \in  {N}_{G}\left( w\right) }\right\}   + \left\{  {w{w}'' :  {w}'' \in  } {V\left( G\right)  \backslash  \left( {{V}_{i}\cup \{ u\} }\right) }\right\}\). Then \({G}^{\prime }\) is non-$r$-partite \({B}_{r,k}\)-free. However,
$$
e\left( {G}^{\prime }\right)  - e\left( G\right)  = n - \left| {V}_{i}\right|  - 1 - {d}_{G}\left( w\right)  \geq  n - \left| {V}_{i}\right|  - 1 - \left( {n - \left| {V}_{i}\right|  - 2}\right)  = 1,
$$
a contradiction to the choice of $G$. 
\end{proof}
Now we come back to show Lemma~\ref{lem3.7}.

Since \(|N_{V_2}(u)\cap N_{V_2}(v)|  \geq  1\), one may take a vertex \(u_2 \in  N_{V_2}(u) \cap N_{V_2}(v)\).
For an integer \(s\) with \(2 \leq  s \leq  r - 1\), assume that for any \(2 \leq  i \leq  s\), there is a 
vertex \({u}_{i}\in V_i\) such that \(\{u,v, u_2, \ldots, u_s\}\) is a clique. 
Consider the number of common neighbors of these vertices in \(V_{s+1}\), by Lemma \ref{lem2.2} and Claims \ref{Claim 4}, \ref{Claim 5}, one has
\begin{align*}
&\,\,\,\,\,\,\,\left|{\left( {{N}_{{V}_{s + 1}}\left( u\right)  \cap  {N}_{{V}_{s + 1}}\left( v\right) }\right)  \bigcap  \left( {\bigcap_{{i = 2}}^{s}{N}_{{V}_{s + 1}}\left( {u}_{i}\right) }\right) }\right|\\
 &\geq
\left| \left( {{N}_{{V}_{s + 1}}\left( u\right)  \cap  {N}_{{V}_{s + 1}}\left( v\right) }\right) \right|  + \sum_{{i = 2}}^{s}\left| {{N}_{{V}_{s + 1}}\left( {u}_{i}\right) }\right|  - \left( {s - 1}\right)  \cdot  \left| {V}_{s + 1}\right|
\\
&\geq  \frac{n}{4r} + \left( {s - 1}\right)  \cdot  \left( {\left| {V}_{s + 1}\right|  - 1}\right)  - \left( {s - 1}\right)  \cdot  \left| {V}_{s + 1}\right|\\
&= \frac{n}{4r} - s + 1 \geq  k.
\end{align*}
That is to say, $u,v, u_2, \ldots, u_s$ have at least $k$ common neighbors in $V_{s+1}$. Thus, for each \(i \in  \left\{2,3,\ldots,r - 1\right\}\), there exists a vertex \({u}_{i} \in  {V}_{i}\) such that \(\left\{  {u,v,{u}_{2},\ldots ,{u}_{r - 1}}\right\}\) is a clique, and they have \(k\) common neighbors \({u}_{r,1},{u}_{r,2},\ldots,{u}_{r,k}\) in \({V}_{r}\). Now \(G\left\lbrack  {\{ u,v,{u}_{2},\ldots ,{u}_{r - 1},{u}_{r,1},{u}_{r,2},\ldots ,{u}_{r,k}\} }\right\rbrack\) is isomorphic to \({B}_{r,k}\), a contradiction. Therefore, there is some $i\in \{2,3,\ldots ,r\}$ such that \({N}_{{V}_{i}}\left( u\right)  \cap  {N}_{{V}_{i}}\left( v\right)  = \emptyset\).
\end{proof}

Now, we are ready to show Theorem \ref{thm1}.
\begin{proof}[\bf Proof of Theorem \ref{thm1}]
By Lemma \ref{lem3.7}, we may assume \({N}_{{V}_{2}}\left( u\right)  \cap  {N}_{{V}_{2}}\left( v\right)  = \emptyset\). Let \({d}_{{V}_{2}}\left( u\right)  = s\), then \({d}_{{V}_{2}}\left( v\right)\)  \(\leq  \left| {V}_{2}\right|  - s\), and \(G\) is a subgraph of \(K_{\left| {V}_{1}\right|,\left| {V}_{2}\right|,\ldots,\left| {V}_{r}\right| }^{s}\), where \(K_{\left| {V}_{1}\right|,\left| {V}_{2}\right|,\ldots,\left| {V}_{r}\right| }^{s}\) is the graph obtained from a complete $r$-partite graph with parts \({V}_{1},{V}_{2},\ldots ,{V}_{r}\) by adding an edge \({uv}\) in \({V}_{1}\) and then deleting \(|V_2| - s\) edges between \(u\) and \({V}_{2}\),\, $s$ edges between \(v\) and \({V}_{2}\) such that \(u\) and \(v\) have no common neighbor in \({V}_{2}\) in the resulting graph. Clearly, \(K_{\left| {V}_{1}\right|,\left| {V}_{2}\right|,\ldots,\left| {V}_{r}\right| }^{s}\) is non-\(r\)-partite \({B}_{r,k}\)-free. Hence, \(e\left( G\right)  \leq  e\left( K_{\left| {V}_{1}\right|,\left| {V}_{2}\right|,\ldots,\left| {V}_{r}\right| }^{s}\right)\), with equality if and only if \(G = K_{\left| {V}_{1}\right|,\left| {V}_{2}\right|,\ldots,\left| {V}_{r}\right| }^{s}\). Note that \(K_{|V_1|, |V_2|,\ldots,|V_r|}^{s} = C_5\left\lbrack  {|V_1|-2,s,1,1}, |V_2|-s\right\rbrack   \vee K_{|V_{3}|,\ldots,|V_r|}\). By Lemma \ref{lem2.4},
$$
e\left( G\right)  \leq  \left({1 - \frac{1}{r}}\right)  \cdot  \frac{{n}^{2}}{2} - \frac{n}{r} + \frac{p\left( {p + 2}\right) }{2r} - \frac{p}{2} + 1
$$
with equality if and only if $G$ is the graph described in Theorem \ref{thm1}(b). 

This completes the proof of Theorem \ref{thm1}. 
\end{proof}

\section{\normalsize Further discussions}%\setcounter{equation}{0}
In this paper, we solve Problem \ref{pb1} for generalized book graphs $B_{r,k},$ where $r\geq 3,\,k\geq1.$ This generalizes a main result in \cite{AFGS,KP2005} in case $n$ is sufficiently large. On the other hand, note that for $r=2$, $B_{2,k}$ is just the book graph $B_k.$ Our result also can be seen as an extension of the main result of Miao, Liu and van Dam \cite{MLD}, which solved Problem \ref{pb1} for $B_k\,(k\geq 2).$ 

Let $G$ be a graph, the \textit{spectral radius} of $G$, denoted by $\rho(G)$, is the largest eigenvalue of the \textit{adjacency matrix} of $G.$ Nikiforov \cite{Nik2009} presented a spectral version of Theorem \ref{thm1.2}: if $H$ is a color-critical graph with $\chi(H)=r+1\,(r\geq 2),$ then $T_r(n)$ is the unique $H$-free graph with maximum spectral radius. Motivated by this, we may pose the following problem, which can be seen as a spectral version of Problem \ref{pb1}.
\begin{pb}\label{pb2}
Given a color-critical graph $H$ with $\chi(H)=r+1$. Determine the graphs with maximum spectral radius among all non-$r$-partite $H$-free graphs.
\end{pb}
Lin, Ning and Wu \cite{LNW} solved Problem \ref{pb2} for $K_3.$ Li and Peng \cite{LP2023} solved Problem \ref{pb2} for $K_{r+1}\,(r\geq 3).$ Recently, Liu and Miao \cite{LM2025} solved Problem \ref{pb2} for $B_{r+1}\,(r\geq 1).$ Motivated by Theorem \ref{thm1}, we will solve Problem \ref{pb2} for $B_{r,k}\,(r\geq 3,k\geq 1)$ in the next step.


\begin{thebibliography}{99}
\small \setlength{\itemsep}{-.8mm}
\bibitem{AFGS} K. Amin, J. Faudree, R.J. Gould, E. Sidorowicz, On the non-\((p-1)\)-partite $K_p$-free graphs, Discuss. Math. Graph Theory 33 (2013) 9-33.
\bibitem{BB1998} B. Bollob\'as, Modern Graph Theory, Springer Science \& Business Media, 1998.  
\bibitem{CJ2002} L. Caccetta, R.-Z. Jia, Edge maximal non-bipartite graphs without odd cycles of prescribed lengths, Graphs Combin. 18 (1) (2002) 75-92.
\bibitem{CFTZ} S. Cioab\u a, L.H. Feng, M. Tait, X.-D. Zhang, The maximum spectral radius of graphs without friendship subgraphs, Electron. J. Combin. 27 (4) (2020) Paper No. 4.22, 19 pp.
\bibitem{E1955} P. Erd\H{o}s, Some theorems on graphs, Riveon Lematematika 9 (1955) 13-17. 
\bibitem{E1966} P. Erd\H{o}s, Some recent results on extremal problems in graph theory (Results), In: Theory of Graphs International Symposium Rome, 1966, Gordon and Breach, New York, Dunod, Paris, 1966, 117-123.
\bibitem{E1968} P. Erd\H{o}s, On some new inequalities concerning extremal properties of graphs, In: Theory of Graphs Proceedings of the Colloquium, Tihany, 1966, Academic Press, New York, 1968, 77-81.
 \bibitem{FS2013} Z. F\"uredi, M. Simonovits, The history of degenerate (bipartite) extremal graph problems, Erd\H{o}s centennial, Bolyai Soc. Math. Stud. 25 (2013) 169-264. 
\bibitem{KP2005} M. Kang, O. Pikhurko, Maximum $K_{r+1}$-free graphs which are not $r$-partite, Mat. Stud. 24 (1) (2005) 12–20.
\bibitem{LP2023} Y.T. Li, Y.J. Peng, Refinement on spectral Tur\'an's theorem, SIAM J. Discrete Math. 37 (4) (2023) 2462-2485.
\bibitem{LNW} H.Q. Lin, B. Ning, B.Y.D.R. Wu, Eigenvalues and triangles in graphs, Combin. Probab. Comput. 30 (2) (2021) 258-270.
\bibitem{LM2025} R.F. Liu, L. Miao, Spectral Tur\'an problem of non-bipartite graphs: forbidden books, European J. Combin. 126 (2025) Paper No. 104136, 18 pp.    
\bibitem{M1907} W. Mantel, Problem 28, Solution by H. Gouwentak, W. Mantel, J. Teixeira de Mattes, F. Schuh and W. A. Wythoff , Wiskundige Opgaven 10 (1907) 60-61. 
\bibitem{MLD} L. Miao, R.F. Liu, E.R. van Dam, Tur\'an number of books in non-bipartite graphs, submitted.       
\bibitem{Nik2009} V. Nikiforov, Spectral saturation: inverting the spectral Tur\'an theorem. Electron. J. Combin. 16 (1) (2009) Research Paper 33, 9 pp.
\bibitem{V2011} V. Nikiforov, Some new results in extremal graph theory: In surveys in Combinatorics 2011, London Math. Society Lecture Note Ser. 392 (2011) 141-181. 
\bibitem{RWWY} S.J. Ren, J. Wang, S.P. Wang, W.H. Yang, A stability result for $C_{2k+1}$-free graphs, SIAM J. Discrete Math. 38 (2) (2024) 1733-1756.
\bibitem{Simo1968} M. Simonovits, A method for solving extremal problems in graph theory, stability problems, Academic Press, New York-London (1968) 279-319. 
\bibitem{Simo1974} M. Simonovits, Extermal graph problems with symmetrical extremal graphs. Additional chromatic conditions, Discrete Math. 7 (1974) 349-376.       
\bibitem{T1941} P. Tur\'an, On an extremal problem in graph theory, Mat. Fiz. Lapok 48 (1941) 436-452 (in Hungarian).     
\bibitem{DW1996} D.B. West, Introduction to Graph Theory, Prentice Hall, Inc., Upper Saddle River, NJ, 1996.         
\end{thebibliography}
\end{document}